\documentclass[11pt]{article}
\usepackage[latin1]{inputenc}
\usepackage{amsmath}
\usepackage{amsfonts}
\usepackage{amssymb}

\newtheorem{theorem}{Theorem}[section]
\def\Rd{\mathbb{R}^d}
\newcommand{\ter}{\hspace{\stretch{3}}$\square$\\[1.8ex]}

\newtheorem{lemma}[theorem]{Lemma}

\newtheorem{prop}[theorem]{Proposition}
\newtheorem{teo}[theorem]{Theorem}
\newtheorem{remark}[theorem]{Remark}

\def\Rd{\mathbb{R}^d}
\def\C{C^+_c(\Rd)}

\def\R{\mathbb{R}}
\def\Ex{\mathbb{E}_x}
\def\E{\mathbb{E}}
\def\S{\mathcal{S}}

\def\L{\mathcal{L}}

\newcommand{\demostracion}{{\it \bf Proof:  }}
\title{\bf Laws of Large Numbers for the Occupation Time of an
Age-Dependent Critical Binary Branching System}
\author{{\sc J. Alfredo L\'{o}pez-Mimbela} \\ {\tt jalfredo@cimat.mx} \\  {\sc Antonio Murillo-Salas}\\ {\tt murillo@cimat.mx} \\
Centro de Investigaci\'{o}n en Matem\'{a}ticas, Guanajuato, Mexico}
\begin{document}
\date{ }
\maketitle
\begin{abstract}
The occupation time of an age-dependent branching particle system in $\Rd$ is
considered, where the initial population is a Poisson random
field and the particles are subject to symmetric $\alpha$-stable
migration, critical binary branching and random lifetimes.
Two regimes of lifetime distributions are considered:
lifetimes with finite mean and lifetimes belonging to the normal
domain of attraction of a $\gamma$-stable law, $\gamma\in(0,1)$. It
is shown that in dimensions $d>\alpha\gamma$ for the
heavy-tailed lifetimes case, and $d>\alpha$ for finite mean lifetimes, the occupation time proccess satisfies a strong law of large
numbers.


\medskip

\end{abstract}

\section{Introduction and background}

In this paper, we obtain laws of large numbers for the occupation
time process of a random population living in the $d$-dimensional
Euclidean space $\Rd$. The evolution of the population is as
follows. Any given individual independently develops a spherically
symmetric $\alpha$-stable process during its lifetime $\tau$, where
$0<\alpha\le2$ and $\tau$ is a random variable having a
non-arithmetic distribution function, and at the end of its life it
either disappears, or  is replaced at the site where it died by two
newborns, each event occurring with probability 1/2.
 The population starts off from a Poisson random field having  Lebesgue measure $\Lambda$ as its intensity. We postulate the usual independence assumptions in branching systems.

 Two regimes
 for the distribution of $\tau$ are considered: either $\tau$
 has  finite mean $\mu>0$, or $\tau$ possesses a
 distribution function $F$
such that
 $F(0)=0$, $F(x)<1$ for all $x\in[0,\infty)$,
and (with $g(u)\sim h(u)$, as $u\rightarrow\infty$, meaning $g(u)/h(u)\rightarrow \mbox{const}$, as $u\rightarrow\infty$)
\begin{equation}\label{tail1}
\bar{F}(u):=1-F(u)\sim u^{-\gamma}{\Gamma(1-\gamma)}^{-1}\quad\mbox{as}\quad
u\longrightarrow\infty,
\end{equation}
where $\gamma\in(0,1)$ and $\Gamma$ denotes the Gamma function,
i.e., $F$ belongs to the normal domain of attraction of a
$\gamma$-stable law. In particular, this allows to consider lifetimes with infinite mean.

Let $X(t)$ denote the simple counting measure on $\Rd$ whose atoms
are the positions of particles alive at  time $t$, and let
 $X\equiv\{X(t)$, $t\ge0\}$. When
$\tau$ has an exponential distribution it is well known that the
measure-valued process $X$ is Markov. In the literature there is a
lot of work about the Markovian model. Our objective here is to
investigate the case when $\tau$ is not necessarily an exponential
random variable, in which case $\{X(t)$, $t\ge0\}$ is no longer a
Markov process. Another striking difference with respect to the case
of exponential lifetimes arises when the particle lifetime
distribution satisfies (\ref{tail1}). When the distribution of
$\tau$ possesses heavy tails, a kind of compensation occurs between
longevity of individuals and clumping of the population:
heavy-tailed lifetimes enhance the mobility of individuals,
favouring in this way the spreading out of particles, and thus
counteracting the clumping of the population. Since clumping goes
along with local extinction (due to critical branching), a smaller
exponent $\gamma$ suits better for stability of the population. As a
matter of fact, Vatutin and Wakolbinger \cite{VW} and Fleischmann,
Vatutin and Wakolbinger \cite{FVW} proved that $X$ admits a
nontrivial equilibrium distribution if and only if $d \ge
\gamma\alpha$. This contrasts with the case of finite-mean (or
exponentially distributed) lifetimes, where the necessary and
sufficient condition for stability is $d>\alpha$. As we shall see,
such qualitative departure from the Markovian model propagates also
to other aspects of the branching particle system, such as the
large-time behavior of its occupation time.

Recall that the {occupation time } of the
measure-valued process $X$ is again a measure-valued process
$J\equiv\{J_t$, $t\geq0\}$, which is  defined by
\begin{equation*}
\langle\psi,J_t\rangle:=\int_0^t\langle\psi,X_s\rangle
ds,\;\;t\geq0,
\end{equation*}
for all bounded measurable function
$\psi:\mathbb{R}^d\rightarrow\mathbb{R}_+$, where the notation $\langle\psi,\nu\rangle$ means $\int\psi\,d\nu$.
Limit theorems for occupation times of  branching systems have been
extensively studied in the context of exponentially distributed lifetimes.
Cox and Griffeath \cite{CG} proved  a strong law of large numbers for the
occupation time of a critical binary branching system.
Moreover, in \cite{CG} it is proved a central limit-type theorem
for the occupation time of the critical binary branching Brownian
motion. M\'{e}l\'{e}ard and Roelly \cite{SS} extended the
law of large numbers of \cite{CG} to branching populations with  general finite-variance critical branching, and quasi-stable particle motions.
Bojdecki, Gorostiza and Talarczyk \cite{BGT,BGT1,BGT2} have investigated the limit fluctuations of
the re-scaled occupation time  $\{J_T(t):=J_{tT},t\geq0\}$ of branching systems,  $T$ being a parameter which tends to infinity. They
have shown that these limits are processes which exhibit
long-range dependence behavior, such as {fractional Brownian
motion} and  {sub-fractional Brownian motion}. See also
\cite{BZ} for related results, where the underlying process is a
branching random walk in the $d$-dimensional lattice.

In this paper we will
prove that, in  dimensions $d>\alpha\gamma$ for
heavy-tailed lifetimes, and $d>\alpha$ for finite-mean lifetimes, the occupation time of the process $X$ satisfies a strong law of
large numbers. Namely, a.s. for any positive continuous function $\psi$ with compact support,
\begin{equation*}
t^{-1}\langle\psi,J_t\rangle\longrightarrow
\langle\psi,\Lambda\rangle \quad \mbox{as
$\quad t\longrightarrow\infty$.}
\end{equation*}
Also, we prove that in  dimensions $d<\alpha\gamma$ for heavy-tailed
lifetimes, and $d<\alpha$ for finite-mean lifetimes, the normalized occupation
time $t^{-1}J_t$ converges to zero a.s. in the sense that, with probability 1, for any ball $A\subset\Rd$ of finite radius,
\begin{equation}\label{null-measure}
t^{-1}\int_0^t1_{\{X_s(A)>0\}}\,ds\longrightarrow 0\quad \mbox{as $\quad
t\longrightarrow\infty$.}
\end{equation}
These results complement ---and partially extend--- those of \cite{CG}
and \cite{SS}. 
We point out that  dimension-dependent behaviors, or parameters,
are  a typical characteristic in this theory because
properties of the branching system treated here are highly related to the
transience-recurrence behavior of the particle motions.
Notice also that, in contrast with the case of
finite-mean lifetimes, in the presence of
heavy-tailed lifetimes the dimension restriction varies according to the decay exponent of the
tail. This phenomenon is reminiscent of the interplay of population clustering and longevity of individuals quoted above.

Our proofs use techniques from \cite{CG}, \cite{Is} and \cite{SS}. To prove the strong law of large numbers in case of heavy-tailed
lifetimes, we first consider the case of ``intermediate dimensions" $\alpha\gamma<d<2\alpha$, and deal afterward
 with ``large dimensions" $d\geq2\alpha$.
Aiming  at applying the Borel-Cantelli lemma, in case of intermediate dimensions we use the re-scaled occupation time process to
upper-bound the variance functional of the occupation time. This
step employs certain  Fourier-transform techniques that we adapted
from \cite{BGT}. We were unable to extend this method to dimensions
$d\ge2\alpha$ due to the lack of proper upper-bounds for the
variance functional of the re-scaled occupation time. To deal with
the case of large dimensions we follow the approach of \cite{SS}. We
use a Markovianized branching system, introduced in Section 4.1
below, which allows us to directly apply the well-known
self-similarity of the symmetric $\alpha$-stable transition
densities. We remark that, in order to use this procedure, we need
to assume that the lifetime distribution possesses a continuous
density function. This contrast with the case of low dimensions,
where no absolute continuity condition is required. We think,
however, that the result should be true for a general lifetime
distributions.

In case of a general non-arithmetic lifetime distribution having
finite mean, our proof of the law of large numbers is carried out
using estimates for the variance functional of the occupation time
process, as well as bounds for the $\alpha$-stable transition
densities. The almost sure convergence
(\ref{null-measure}) is proved by combining Borel-Cantelli's
Lemma with some estimates from \cite{VW} related to extinction
probabilities.

The analysis at the ``critical dimensions" $d=\alpha\gamma$ in the
heavy-tailed case, and $d=\alpha$ for finite mean lifetimes, is much
more difficult to carry out, as can be seen from \cite{FG}, where
the occupation time (at the critical dimension $d=\alpha/\beta$) of
the so-called $(d,\alpha,\beta)$-superprocess is considered. The
approach there strongly relies on the classical semilinear equation
characterizing the Laplace functional of the occupation time, see
Lemma 3.4 in \cite{FG}. In our case, due to the non-exponential
lifetimes, we do not have the above-mentioned equation. Thus, laws
of large numbers for our model at critical dimensions remain to be
investigated.
\section{Laws of large numbers}

Following \cite{BGT}, we define the re-scaled occupation time
process $\{J_T(t):=J_{tT},\,t\geq0\}$, i.e., for any positive
bounded measurable function $\varphi$,
\begin{equation}\label{rescaled}
\langle\varphi,J_T(t)\rangle:=\int_0^{tT}\langle\varphi,X_s\rangle
ds=T\int_0^t\langle\varphi,X_{sT}\rangle ds,\;\;t\geq0.
\end{equation}
Notice that, by Fubini's theorem,
\begin{equation}\label{meanR-OccT}
\E\langle\varphi,J_T(1)\rangle=\langle\varphi,\Lambda\rangle T,
\end{equation}
since $\E\langle\varphi,X_t\rangle=\langle\varphi,\Lambda\rangle$.
We remark that studying the asymptotic behavior of $\langle
\varphi,J_t\rangle/t$ as $t\longrightarrow\infty$, is the same
 as investigating the asymptotic behavior of $\langle
\varphi,J_T(1)\rangle/T$ as $T\longrightarrow\infty$.

In what follows, $\C$ denotes the space of nonnegative continuous
functions $\varphi:\Rd\longrightarrow\R_+$ with compact support. The
main results of this paper are the following theorems.

\begin{teo}\label{lowdimensionsT} Let $F$ be a non-arithmetic
distribution function satisfying (\ref{tail1}).
\\(a) Assume that $\alpha\gamma<d<2\alpha$. Then, a.s. for any
$\varphi\in\C$,
\begin{equation}
T^{-1}\langle\varphi,J_T(1)\rangle\longrightarrow\langle\varphi,\Lambda\rangle\quad\mbox{as}\quad
T\longrightarrow\infty.
\end{equation}
(b) Suppose that $d\geq2\alpha$ and $F$ possesses a continuous density $f$.
Then, a.s. for all $\varphi\in\C$,
\begin{equation}
T^{-1}\langle\varphi,J_T(1)\rangle\longrightarrow\langle\varphi,\Lambda\rangle\quad\mbox{as}\quad
T\longrightarrow\infty.
\end{equation}
\end{teo}
\noindent Our next theorem complements the law of large numbers of
\cite{CG} and \cite{SS}, which were proved only in the case of
exponentially distributed lifetimes.

\begin{teo}\label{LLNFinitemean} Assume that $d>\alpha$, and let $F$
be a non-arithmetic distribution function with finite mean $\mu>0$.
Then, a.s. for any $\varphi\in\C$,
\begin{equation}
T^{-1}\langle\varphi,J_T(1)\rangle\longrightarrow\langle\varphi,\Lambda\rangle\quad\mbox{as}\quad T\longrightarrow\infty.
\end{equation}
\end{teo}

\begin{teo}\label{subcritical}
Let $F$ be a non-arithmetic
distribution function which satisfies (\ref{tail1}), and assume that $d<\alpha\gamma$. Then, a.s. for any ball $A\subset\Rd$ of finite radius,
\begin{equation}\label{subcritical1}
T^{-1}\int_0^T1_{\{X_s(A)>0\}}\,ds\longrightarrow0\quad\mbox{as}\quad T\longrightarrow\infty.
\end{equation}
\end{teo}
\begin{remark} (a) Notice that condition $\alpha\gamma<d<2\alpha$ in Theorem \ref{lowdimensionsT}
allows $d\leq\alpha$, which contrasts with the classical case of
exponentially distributed lifetimes, where $d>\alpha$.
\\
(b) When the particle lifetimes have an exponential distribution
with mean $\lambda^{-1}$, Theorem \ref{LLNFinitemean} reduces to
Theorem 4 of \cite{SS}.
\\(c) In case of low dimensions $d<\alpha\gamma$, a genuine counterpart to Theorem \ref{lowdimensionsT}  would be a statement ensuring a.s. vague convergence of
$T^{-1}J_T(1)$ to the zero measure as $T\to\infty$. This was proved by Sawyer and Fleischman \cite{S-F} for the occupation time of
critical branching Brownian motion (see also \cite{Is} for a related result regarding super-Brownian motion's occupation time). For our model, here we
prove only the slightly weaker result (\ref{subcritical1}), which implies that, with probability 1, the proportion of time that the branching system charges any given bounded set  vanishes asymptotically as $T\to\infty$. The specific form (\ref{subcritical1}) of a.s. convergence  was suggested to us by an anonymous referee.
\\ (d) The extent of dimensions in our results is narrow due to our
choice of critical, binary-branching mechanism. A less restrictive
assumption, such as critical $(1+\beta)$-branching, $\beta\in(0,1)$,
would expand the dimension range.

\end{remark}

\section{Some moment calculations}

Let  $Z_t(A)$ denote the number of individuals living in
$A\in\mathcal{B}(\Rd)$ at time $t$, in a population starting with
one particle at time $t=0$. Following \cite{KS} we define
\begin{equation}\label{Iniciauno}
Q_t\varphi(x):=\mathbb{E}_x\left[1-e^{-\langle\varphi,Z_t\rangle}\right],\quad
x\in\mathbb{R}^d,\quad  t\geq0,
\end{equation}
where  $\varphi\in C^+_c(\Rd)$ and $\E_x$ means that the initial particle is located at $x\in\Rd$. Since the initial population $X_0$ is
Poissonian, we have
\begin{eqnarray}\label{uno}
\mathbb{E}e^{-\langle\varphi,
X_t\rangle}&=&\nonumber\exp\left(-\int\mathbb{E}_x\left[1-e^{-\langle\varphi,
Z_t\rangle}\right]dx\right)\\&=&\exp\left(-\int
Q_t\varphi(x)dx\right),\;\;\varphi\in\C.
\end{eqnarray}
Let $\{\tau_k,k\geq1\}$
be a sequence of i.i.d. random variables with common distribution function
$F$, and let
$$
N_t=\sum_{k=1}^\infty1_{\{S_k\leq t\}},\quad t\geq0,
$$
where the random sequence $\{S_k,k\geq0\}$ is recursively defined by
$$
S_0=0,\quad S_{k+1}=S_k+\tau_k,\quad k\geq0.
$$
For any $p=1,2,\ldots$, $0<t_p\leq t_{p-1},\ldots,t_1<\infty$,
$\varphi_1,\varphi_2,\ldots,\varphi_p\in \C$ and
$\theta_1,\ldots,\theta_p\in\R$, we define
$\bar{t}=(t_1,t_2,\ldots,t_p)$,
$\bar{t}-s=(t_1-s,t_2-s,\ldots,t_p-s)$,
$\theta_{(p)}=(\theta_1,\ldots,\theta_p)^\prime$ and
\begin{equation*}
Q^p_{\bar{t}}\theta_{(p)}(x)=\mathbb{E}_x\left[1-e^{-\sum_{j=1}^p\theta_j\langle\varphi_j,Z_{t_j}\rangle}\right].
\end{equation*}
Let $\{B_s$, $s\ge0\}$ denote the spherically symmetric
$\alpha$-stable process in $\Rd$, with transition density functions
$\{p_t(x,y)$, $t>0$, $x,y\in\mathbb{R}^d\}$, and semigroup
$\{\mathcal{S}_t$, $t\geq0\}$.
Our moment calculations use the following result which is borrowed from \cite{KS} (Section 4.3), and which we include for convenience.

\begin{prop}\label{finite-dimensional} The function $Q^p_{\bar{t}}\theta_{(p)}$ satisfies
\begin{eqnarray*}
Q^p_{\bar{t}}\theta_{(p)}(x)&=&\mathbb{E}_x\left[1-e^{-\sum_{j=1}^p\theta_j\varphi_j(B_{t_j})}-\int_0^{t_p}\frac{1}{2}\left(Q^p_{\bar{t}-s}\theta_{(p)}(B_s)\right)^2dN_s\right.\\&
&-\left.\sum_{i=1}^{p-1}\left(1-e^{-\sum_{j=i+1}^p\theta_j\varphi_j(B_{t_j})}\right)\int_{t_{i+1}}^{t_i}\frac{1}{2}\left(Q^i_{\bar{t}-s}\theta_{(i)}(B_s)\right)^2dN_s\right].
\end{eqnarray*}
\end {prop}
As in (\ref{uno}), since the initial population is
Poissonian we have
\begin{eqnarray}\label{finite-dimesional1}
\mathbb{E}\left[e^{-\sum_{j=1}^p\theta_j\langle\varphi_j,X_{t_j}\rangle}\right]&=&\exp\left(-\int\mathbb{E}_x\left[1-e^{-\sum_{j=1}^p\theta_j\langle\varphi,Z_{t_j}\rangle}\right]dx\right)\nonumber\\&=&\exp\left(-\int
Q^p_{\bar{t}}\theta_{(p)}(x)dx\right).
\end{eqnarray}
Using criticality of the branching, and that Lebesgue measure is invariant for the semigroup of the symmetric $\alpha$-stable process,  it is easy to see that
\begin{equation}\label{media}
m(t,\varphi):=\E[\langle\varphi,X_t\rangle]=
\langle\varphi,\Lambda\rangle,\quad t\geq0,\quad\varphi\in\C.
\end{equation}

\begin{lemma}\label{COVP} Let $0<s\leq t<\infty$ and $\psi,\varphi\in \C$. Then,
\begin{eqnarray}
\nonumber
C_x(s,\varphi;t,\psi)&:=&\Ex\left[\langle\varphi,Z_{s}\rangle\langle\psi,Z_{t}\rangle\right]
\\&=&\mathbb{E}_x\left[\varphi(B_s)\psi(B_t)+\int_0^sm_{B_r}(t-r,\psi)m_{B_r}(s-r,\varphi)dN_r\right],
\end{eqnarray}
where $m_x(t,\varphi)=\Ex[\langle\varphi,Z_t\rangle]$.
\end{lemma}
\demostracion In order to preserve the notation in Proposition
\ref{finite-dimensional}, we put $p=2$, $t_1=t$, $t_2=s$,
$\varphi_1=\psi$ and $\varphi_2=\varphi$. Then we have
$$
C_x(t_1,\varphi_1;t_2,\varphi_2)=-\frac{\partial^2}{\partial\theta_1\partial\theta_2}Q^2_{\bar{t}}\theta_{(2)}(x)\bigg|_{\theta_1=\theta_2=0^+},
$$
where
\begin{eqnarray*}
\frac{\partial^2}{\partial\theta_1\partial\theta_2}Q^2_{\bar{t}}\theta_{(2)}(x)&=&\mathbb{E}_x\bigg[-\varphi_1(B_{t_1})\varphi_2(B_{t_2})e^{-\theta_1\varphi(B_{t_1})-\theta_2\varphi_2(B_{t_2})}\\&
&-\int_0^{t_2}\frac{\partial}{\partial\theta_2}Q^2_{\bar{t}-r}\theta_{(2)}(B_r)\frac{\partial}{\partial\theta_1}Q^2_{\bar{t}-r}\theta_{(2)}(B_r)dN_r\\&
&-\int_0^{t_2}\left(Q^2_{\bar{t}-r}\theta_{(2)}(B_r)\right)\frac{\partial^2}{\partial\theta_2\partial\theta_1}Q^2_{\bar{t}-r}\theta_{(2)}(B_r)dN_r\\&
&\left.-\varphi_2(B_{t_2})e^{-\theta\varphi_2(B_{t_2})}\int_{t_1}^{t_2}\left(Q^1_{t_2-r}\theta_1(B_r)\right)\frac{\partial}{\partial\theta_1}Q^1_{t_2-r}\theta_1(B_r)dN_r\right].
\end{eqnarray*}
Evaluating at $\theta_1=\theta_2=0$  we finish the proof.\hfill$\Box$

\begin{prop} \label{COVp} Let $0<s\leq t<\infty$ and $\psi,\varphi\in\C$. Then,
\begin{equation}\label{COVp-1}
C(s,\varphi;t,\psi)
:=\mbox{\rm Cov}\left(\langle\varphi,X_s\rangle,\langle\psi,X_t\rangle\right)=
\langle\varphi\S_{t-s}\psi,\Lambda\rangle+\int_0^s\langle\left(\S_{s-r}\varphi\right)\left(\S_{t-r}\psi\right),\Lambda\rangle
\,dU(r),
\end{equation}
where $U(r)=\sum_{k=0}^\infty F^{*k}(r)$.
\end{prop}
\demostracion We put $p=2$ in
(\ref{finite-dimesional1}) and use the same notations as in the proof of Lemma \ref{COVP}. Then,
\begin{eqnarray*}
\mathbb{E}\left[\langle\varphi_1,X_{t_1}\rangle\langle\varphi_2,X_{t_2}\rangle\right]&=&
\frac{\partial^2}{\partial\theta_1\partial\theta_2}\exp\left(-\int
Q^2_{\bar{t}}\theta_{(2)}(x)\,dx\right)\bigg|_{\theta_1=\theta_2=0^+}\\
&=&\left[-\frac{\partial^2}{\partial\theta_1\partial\theta_2}\int
Q^2_{\bar{t}}\theta_{(2)}(x)dx\right.\\& &\left.+\int
\frac{\partial}{\partial\theta_1}
Q^2_{\bar{t}}\theta_{(2)}(x)\,dx\int
\frac{\partial}{\partial\theta_2}
Q^2_{\bar{t}}\theta_{(2)}(x)\,dx\right]\bigg|_{\theta_1=\theta_2=0^+}\\
&=&\int
C_x(t_1,\varphi_1;t_2,\varphi_2)dx+\int m_x(t_1,\varphi_1)dx\int
m_x(t_2,\varphi_2)dx,
\end{eqnarray*}
and from Lemma \ref{COVP} we obtain
\begin{equation}\label{COV}
C(s,\varphi;t,\psi)=\int_{\mathbb{R}^d}
\mathbb{E}_x\left[\varphi(B_s)\psi(B_t)+\int_0^sm_{B_r}(t-r,\psi)m_{B_r}(s-r,\varphi)dN_r\right]dx,
\end{equation}
which completes the proof.\hfill$\Box$

\section{Markovianizing an age-dependent branching system}\label{SEC4}
In this section we introduce a Markovian  branching system
$\{\bar{X}_t,\;t\geq0\}$ which will be used to prove Theorem
\ref{lowdimensionsT} (b). Let $X\equiv\{X_t$, $t\ge0\}$ be the
branching system defined in Section 1. For any $t\ge0$, let
$\bar{X}_t$ denote the population in $\R_+\times\Rd$ ($\R_+=[0,\infty)$) obtained by
attaching to each individual $\delta_x\in X_t$ its age. Namely, for
each $t\geq0$,
\begin{equation}\label{asterisco1}
\bar{X}_t=\sum_i\delta_{(\eta^i_t,\xi^i_t)},
\end{equation}
where $\eta^i_t$ and $\xi^i_t$ denotes respectively, the age and
position of the $i^{th}$ particle at time $t$, and the summation is
over all particles alive at time $t$.  Let us assume that
$\bar{X}_0$ is a Poisson random field on $\R_+\times\Rd$ with
intensity measure $F\times\Lambda$. Here, $F$ also means the
Lebesgue-Stieltjes measure corresponding to $F$. The probability
generating function of the branching law is denoted by $\Phi$. Thus,
for critical binary branching, $\Phi(s)\equiv\frac{1}{2}(1+s^2)$,
$-1\leq s\leq1$.

Given a counting measure $\nu$ on $\R_+\times\Rd$,  and a measurable function
$\phi:\R_+\times\Rd\longrightarrow(0,1]$,  we
define
\begin{equation*}
G_\phi(\nu):=\exp\left(\langle\log\phi,\nu\rangle\right).
\end{equation*}
It can be shown that the infinitesimal generator of
$\{\bar{X}_t,\;t\geq0\}$ evaluated at the function $G_\phi(\nu)$ is
given by
\begin{equation}\label{GI-BPS}
\mathcal{G}G_\phi(\nu)=G_\phi(\nu)\left\langle
\frac{\mathcal{L}\phi(*,\cdot)+\lambda(*)[\Phi(\phi(0,\cdot))
-\phi(0,\cdot)]}{\phi(*,\cdot)},\nu\right\rangle,
\end{equation}
where
\begin{equation}
\label{asterisco2}
\lambda(u)=\frac{f(u)}{1-F(u)},\;\;u\geq0,
\end{equation}
is the {hazard rate function} associated to $F$,
and \begin{equation}\label{B-IG}
\mathcal{L}\phi(u,x)=\frac{\partial\phi(u,x)}{\partial
u}+\Delta_\alpha\phi(u,x)-\lambda(u)\left[\phi(u,x)-\phi(0,x)\right],
\end{equation}
where the function $\phi$ is such that $\phi(\cdot,x)\in
C^1_b(\R_+)$ for any $x\in\Rd$, and $\phi(u,\cdot)\in
C^\infty_c(\Rd)$ for any $u\in\R_+$. Here $C^1_b(\R_+)$ denotes the
set of all bounded functions with continuous first derivative, and
$C^\infty_c(\Rd)$ denotes the space of infinitely differentiable
functions from $\Rd$ to $\R$, having compact support.  The operator
$\L$ is the infinitesimal generator of a Markov process on
$\R_+\times\Rd$ whose semigroup is denoted by $\{\tilde{T}_t,t\geq0\}$,
see \cite{M} for details.

\begin{prop}
\label{LTOC}Let $\bar{X}\equiv\{\bar{X}_t,\,t\geq0\}$ as before and
let $\bar{X}_0$ be a Poisson random field on $\R_+\times\Rd$ with
intensity measure $F\times\Lambda$. The joint Laplace functional of
the branching particle system $\bar{X}$ and its occupation time is
given by
\begin{equation*}
\mathbb{E}\left[e^{-\langle\psi,\bar{X}_t\rangle-\int_0^t\langle\phi,\bar{X}_s\rangle
ds}\right] = e^{-\langle V_t^\psi\phi,
F\times\Lambda\rangle},\;\;t\geq0,
\end{equation*}
for all measurable functions
$\psi,\phi:\R_+\times\Rd\longrightarrow\R_+$ with compact support,
where $V_t^\psi\phi$ satisfies, in the mild sense, the non-linear
evolution equation
\begin{eqnarray}
\nonumber\frac{\partial}{\partial t}V_t^\psi\phi(u,x) &=&
\mathcal{L}V_t^\psi\phi(u,x)-\lambda(u)[\Phi(1-V_t^\psi\phi(0,x))-(1-V_t^\psi\phi(0,x))]\\&
&\label{TOC}+\phi(u,x)(1-V_t^\psi\phi(u,x)),\\\nonumber
V_0^\psi\phi(u,x) &=& 1-e^{-\psi(u,x)}.
\end{eqnarray}
\end{prop}
\demostracion The proof is carried out using the martingale problem
for $\{\bar{X}_t,\,t\geq0\}$, and It\^{o}'s formula. We omit the details.

\section{Proof of Theorem \ref{lowdimensionsT}}
 We shall prove
the law of large numbers in two steps. First we show that the
result holds for intermediate dimensions $\alpha\gamma<d<2\alpha$; this
part of the proof relies on the non-Markovian branching system defined
in Section 1, and uses upper bounds for the covariance
functional. In the second step, we consider ``large" dimensions
$d\geq2\alpha$, and in this case we use the Markovianized branching
system described in Section \ref{SEC4}.

\subsection{Proof of Theorem \ref{lowdimensionsT} (a)}
In this section we assume that $\alpha\gamma<d<2\alpha$.
\begin{lemma} \label{cota}Suppose that the hypothesis in Theorem \ref{lowdimensionsT} hold. Then, for each $\epsilon>0$ and all $T>0$ large
enough,
\begin{equation*}
P\left(\left|T^{-1}\langle\varphi,J_T(1)\rangle-\langle\varphi,
\Lambda\rangle\right|>\epsilon\right)\leq\frac{2}{\epsilon^2}\left(
c_3T^{-2}+c_1T^{-1}+c_2T^{-d/\alpha}+c_4 T^{\gamma-d/\alpha}\right),
\end{equation*}
for some positive constants $c_1,\ldots,c_4$.
\end{lemma}
\demostracion Let $\epsilon>0$ be given. Then, using Chebyshev's
inequality and (\ref{meanR-OccT}),
\begin{eqnarray*}
P\left(\left|T^{-1}\langle\varphi,J_T(1)\rangle-\langle\varphi,
\Lambda\rangle\right|>\epsilon\right)&
\leq
&\frac{1}{\epsilon^2}\E\left(T^{-1}\langle\varphi,J_T(1)\rangle-\langle\varphi,\Lambda\rangle\right)^2\\&
=
&\frac{1}{\epsilon^2T^2}\mbox{Cov}\left(\langle\varphi,J_T(1)\rangle,\langle\varphi,J_T(1)\rangle\right)\\&=&\frac{1}{
\epsilon^2}\int_0^1\int_0^1\mbox{Cov}\left(\langle\varphi,X_{sT}\rangle,\langle\varphi,X_{tT}\rangle\right)dt\,ds,
\end{eqnarray*}
where the last equality follows from (\ref{rescaled}). By changing the order of integration we obtain that
\begin{equation}
P\left(\left|T^{-1}\langle\varphi,J_T(1)\rangle-
\langle\varphi,\Lambda\rangle\right|>\epsilon\right)\leq\frac{2}{\epsilon^2}\int_0^1dv\int_0^v\mbox{Cov}\left(\langle\varphi,X_{uT}\rangle,\langle\varphi,X_{vT}\rangle\right)du.
\end{equation}
Therefore, from Proposition \ref{COVp} we deduce that
\begin{equation}
P\left(\left|T^{-1}\langle\varphi,J_T(1)\rangle-
\langle\varphi,\Lambda\rangle\right|>\epsilon\right)\ \leq \
(I)+(II),
\end{equation}
with
\begin{equation*}
(I):=\frac{2}{\epsilon^2}\int_0^1dv\int_0^vdu\langle\varphi
\S_{T(v-u)}\varphi,\Lambda\rangle
\end{equation*}
and
\begin{equation*}
(II):=\frac{2}{\epsilon^2}\int_0^1dv\int_0^vdu\int_0^udU(Tr)\langle\varphi\S_{T(v+u-2r)}\varphi,\Lambda\rangle,
\end{equation*}
where, to obtain (II), we used self-adjointness of $\S_t$ with respect to $\Lambda$, $t\ge0$.
Our next goal is to derive useful upper bounds for the two integrals (I)
and (II). Firstly, by performing the change of variables $s=(v-u)T$
and $t=vT$, we get that
\begin{eqnarray*}
\frac{\epsilon^2}{2}(I)&=&\frac{1}{T^2}\int_0^Tdt\int_0^tds\langle\varphi\S_s\varphi,
\Lambda\rangle\\&=&\frac{1}{T^2}\int_0^Adt\int_0^tds\langle\varphi\S_s\varphi,\Lambda\rangle+\frac{1}{T^2}
\int_A^Tdt\int_0^tds\langle\varphi\S_s\varphi,\Lambda\rangle
\end{eqnarray*}
for any $A>0$, where
\begin{eqnarray*}
\int_0^t\langle\varphi\S_s\varphi,\Lambda\rangle
ds&=&\int_0^t\int_{\Rd}\int_{\Rd}\varphi(x)p_s(x-y)\varphi(y)dy\,dx\,ds\\
&=&\int_{\Rd}\int_{\Rd}\varphi(x)\varphi(y)\int_0^tp_s(x-y)ds\,dy\,dx\\
&\leq&\int_{\Rd}\int_{\Rd}\varphi(x)\varphi(y)\mbox{const.}\left(|x-y|^{\alpha-d}+t^{1-d/\alpha}\right)dy\,dx
\end{eqnarray*}
since
\begin{equation*}
\int_0^tp_s(x-y)ds\leq
\mbox{const.}\left(|x-y|^{\alpha-d}+t^{1-d/\alpha}\right)
\end{equation*}
 due to self-similarity of the $\alpha$-stable transition densities. Notice that
\begin{equation*}
\int_{\Rd}\int_{\Rd}\varphi(x)\varphi(y)|x-y|^{\alpha-d}dy\,dx<\infty,
\end{equation*}
which, for $d>\alpha$, follows from Lemma 5.3 in \cite{Is}. Hence,
\begin{eqnarray*}
\int_A^Tdt\int_0^t\langle\varphi\S_s\varphi,\Lambda\rangle
ds&\leq&\mbox{const.}\int_{\Rd}\int_{\Rd}\varphi(x)\varphi(y)\int_A^T\left(|x-y|^{\alpha-d}+t^{1-d/\alpha}
\right)dt\,dy\,dx\\&=&c_1(T-A)+c_2(T^{1-d/\alpha}-A^{1-d/\alpha})
\end{eqnarray*}
for some constants $c_1,c_2>0$. Therefore,
\begin{equation}\label{cotaI}
(I)\leq\frac{2}{\epsilon^2}\left(\frac{c_3}{T^2}+c_1\frac{T}{T^2}+c_2\frac{T^{2-d/\alpha}}{T^2}\right),
\end{equation}
where
$$c_3=\int_0^Adt\int_0^t\langle\varphi\S_s\varphi,\Lambda\rangle ds.$$

Before estimating the integral $(II)$, we recall that $U(t)\sim
t^{\gamma}/\Gamma(1+\gamma)$ as $t\rightarrow\infty$ because of
$\bar{F}(t)\sim t^{-\gamma}/\Gamma(1-\gamma)$, see \cite{B.et.al.},
p. 361. Then, writing $\hat{\varphi}$ for the Fourier transform of
$\varphi$, we obtain
\begin{eqnarray*}
\frac{\epsilon^2}{2}(II)&=&\int_0^1dv\int_0^vdu\int_0^udU(Tr)\frac{1}{(2\pi)^d}\int_{\Rd}
dy|\hat{\varphi}(y)|^2e^{-T(v+u-2r)|y|^\alpha}\\&=&\frac{1}{(2\pi)^d}\int_0^1dv\int_0^vdu\int_{\Rd}
dy|\hat{\varphi}(y)|^2\int_0^udU(Tr)e^{-T(v+u-2r)|y|^\alpha}\\&\sim&\frac{\gamma
T^\gamma}{\Gamma(1+\gamma)(2\pi)^d}\int_0^1dv\int_0^vdu\int_{\Rd}
dy|\hat{\varphi}(y)|^2\int_0^ue^{-T(v+u-2r)|y|^\alpha}r^{\gamma-1}dr,
\end{eqnarray*}
and, after the change of variables $z=(T(v+u-2r))^{1/\alpha}y$, we
conclude that
\begin{eqnarray*}
\frac{\epsilon^2}{2}(II)&\sim&\frac{\gamma
T^\gamma}{\Gamma(1+\gamma)(2\pi)^d}\int_0^1dv\int_0^vdu\int_{\Rd}
dz\int_0^uT^{-d/\alpha}(v+u-2r)^{-d/\alpha}\\&
&\times|\hat{\varphi}(T^{-d/\alpha}(v+u-2r)^{-d/\alpha}z)|^2e^{-|z|^\alpha}r^{\gamma-1}dr\\&\leq&\frac{\gamma
T^{\gamma-d/\alpha}\langle\varphi,\Lambda\rangle^2}{\Gamma(1+\gamma)(2\pi)^d}\int_{\Rd}dze^{-|z|^\alpha
}\int_0^1dv\int_0^vdu\int_0^u(v+u-2r)^{-d/\alpha}r^{\gamma-1}dr,
\end{eqnarray*}
 where to obtain the last inequality we have used the well known fact that $|\hat{\varphi}(z)|\leq(2\pi)^{-d}\langle|\varphi|,\Lambda\rangle$ for any $L^1$-function $\varphi$.
Changing the order of integration into the above expression yields
\begin{eqnarray*}
&=&\frac{\gamma
T^{\gamma-d/\alpha}\langle\varphi,\Lambda\rangle^2}{\Gamma(1+\gamma)(2\pi)^d}\int_{\Rd}dze^{-|z|^\alpha}\int_0^1dv\int_0^vr^{\gamma-1}\int_r^vdu(u+v-2r)^{-d/\alpha}dr\\&=&\frac{\gamma
T^{\gamma-d/\alpha}\langle\varphi,\Lambda\rangle^2}{\Gamma(1+\gamma)(2\pi)^d}\int_{\Rd}dze^{-|z|^\alpha
}\int_0^1dv\int_0^vr^{\gamma-1}\frac{2^{1-d/\alpha}(v-r)^{1-d/\alpha}-(v-r)^{1-d/\alpha}}{1-d/\alpha}dr\\&=&\frac{\gamma
T^{\gamma-d/\alpha}\langle\varphi,\Lambda\rangle^2}{\Gamma(1+\gamma)(2\pi)^d}\frac{2^{1-d/\alpha}-1}{1-d/\alpha}\int_{\Rd}dze^{-|z|^\alpha
}\int_0^1dv\int_0^vr^{\gamma-1}(v-r)^{1-d/\alpha}\,dr.
\end{eqnarray*}
 Notice that in the above calculations we have implicitly assumed that $d\ne\alpha$. The case $d=\alpha$ can be treated in a similar way (and renders the same conclusion  (\ref{ast})). Changing again the order of integration we get
\begin{eqnarray*}
& =&\frac{\gamma
T^{\gamma-d/\alpha}\langle\varphi,\Lambda\rangle^2}{\Gamma(1+\gamma)(2\pi)^d}\frac{2^{1-d/\alpha}-1}{1-d/\alpha}\int_{\Rd}dze^{-|z|^\alpha
}\int_0^1r^{\gamma-1}\int_r^1dv(v-r)^{1-d/\alpha}dr\\&=&\frac{\gamma
T^{\gamma-d/\alpha}\langle\varphi,\Lambda\rangle^2}{\Gamma(1+\gamma)(2\pi)^d}\frac{2^{1-d/\alpha}-1}{(1-d/\alpha)(2-d/\alpha)}
\int_{\Rd}dze^{-|z|^\alpha
}\int_0^1r^{\gamma-1}(1-r)^{2-d/\alpha}dr,
\end{eqnarray*}
where the last equality is finite since by assumption $d<2\alpha$.
Hence, for $T$ large enough
\begin{equation}\label{ast}
(II)\leq \frac{2}{\epsilon^2}c_4 T^{\gamma-d/\alpha}.
\end{equation}
Therefore,
\begin{equation*}
P\left(|T^{-1}\langle\varphi,J_T(1)\rangle-\langle\varphi,\Lambda\rangle|>\epsilon\right)\leq
\frac{2}{\epsilon^2}\left( c_3T^{-2}+c_1T^{-1}+c_2T^{-d/\alpha}+c_4
T^{\gamma-d/\alpha}\right).
\end{equation*}
\ter

\noindent\textbf{Proof of Theorem \ref{lowdimensionsT} (a):} Let
$\epsilon>0$ and $a>1$ be given constants, and let $T_n=a^n$ for
$n=1,2,\dots$. Then,
\begin{eqnarray*}
\lefteqn{\sum_{n=1}^\infty
P\left(|T_n^{-1}\langle\varphi,J_{T_n}(1)\rangle-\langle\varphi,\Lambda\rangle|>\epsilon\right)}\\&\leq&\frac{2}{\epsilon^2}\sum_{n=1}^\infty
\left(c_3T_n^{-2}+c_1T_n^{-1}+c_2T_n^{-d/\alpha}+c_4T_n^{\gamma-d/\alpha}\right)<\infty
\end{eqnarray*}
due to the assumption $d>\gamma\alpha$. Therefore, for any given
$\varphi\in\C$, a.s.,
\begin{equation*}
T_n^{-1}\langle\varphi,J_{T_n}(1)\rangle\longrightarrow\langle\varphi,\Lambda\rangle\quad\mbox{as}\quad n\longrightarrow\infty.
\end{equation*}
Now we observe that, for each $T>1$, there exists some non-negative
integer $n(T)$ such that $a^{n(T)}\leq T\leq a^{n(T)+1}$, and that
$n(T)\longrightarrow\infty$ as $T\longrightarrow\infty$. Hence,
\begin{equation*}
\frac{\langle\varphi,J_{a^{n(T)}}(1)\rangle}{a^{n(T)+1}}\leq\frac{\langle\varphi,J_{T}(1)\rangle}{T}
\leq\frac{\langle\varphi,J_{a^{n(T)+1}}(1)\rangle}{a^{n(T)}},
\end{equation*}
and
\begin{equation*}
\frac{\langle\varphi,\Lambda\rangle}{a}\leq\liminf_{T\rightarrow\infty}\frac{\langle\varphi,J_{T}(1)
\rangle}{T}\leq\limsup_{T\rightarrow\infty}\frac{\langle\varphi,J_{T}(1)\rangle}{T}\leq\langle\varphi,\Lambda\rangle
a,
\end{equation*}
these inequalities being true for any $a>1$. Letting $a\rightarrow1$
yields that
$$\lim_{T\longrightarrow\infty}T^{-1}\langle\varphi,J_t(1)\rangle=\langle\varphi,\Lambda\rangle$$
a.s., where the null set (where the limit may not exist) depends
upon $\varphi$. Nonetheless, a null set can be chosen not to depend
on $\varphi$ as is the proof of Theorem 1 in \cite{Is1}.
\hfill$\Box$

\subsection{Proof of Theorem \ref{lowdimensionsT} (b)}
Throughout this section we assume that $d\geq2\alpha$.
 The proof of part (b) in Theorem \ref{lowdimensionsT}   follows, as
in part (a), from  Chebyshev's inequality
\begin{equation}\label{c1}
P\left\{\frac{|\langle\phi,J_t\rangle -
\langle\phi,\Lambda\rangle|}{t}>\epsilon\right\}
\leq\frac{1}{t^2\epsilon^2}\mbox{Var}\langle\phi,J_t\rangle,\quad
t\geq0,\quad\epsilon>0
\end{equation}
and  Lemma \ref{MEAN} below. Recall that $\lambda$ is defined in (\ref{asterisco2}).

\begin{lemma}\label{MEAN} i) Let $\phi:\mathbb{R}^d\rightarrow\mathbb{R}_+$ be
a measurable function with compact support. Then, for each $t\geq0$,
\begin{equation}\label{asterisco3}
\mathbb{E}\langle\phi,J_t\rangle=\langle\phi,\Lambda\rangle t,
\end{equation}
and
\begin{equation}\label{CotaVar}
\emph{\mbox{Var}}\langle\phi,J_t\rangle \leq
\langle\lambda,F\rangle\emph{\mbox{Const}}(\phi)(t +
t^{3-d/\alpha}) + 2\emph{\mbox{Const}}(\phi)(t +
t^{2-d/\alpha}).
\end{equation}
\end{lemma}
\demostracion First we prove (\ref{asterisco3}). For the given function $\phi$
we define the extended function $\bar{\phi}(u,x)\equiv\phi(x)$,
$(u,x)\in\R_+\times\Rd$. Then, for any $k\geq0$ we define
\begin{equation}\label{LFTOC}
L_t(k\bar{\phi}) =
\mathbb{E}\left[e^{-k\int_0^t\langle\bar{\phi},\bar{X}_s\rangle
ds}\right]\nonumber =  e^{-\langle
V_t(k\bar{\phi}),F\times\Lambda\rangle},
\end{equation}
where $V_t(k\bar{\phi})$ satisfies (\ref{TOC}) with $\bar{\phi}$
substituted by $k\bar{\phi}$, and $\psi\equiv0$. Notice that
\begin{eqnarray*}
\mathbb{E} \langle{\phi},J_t\rangle& = &
-\frac{d}{dk}\mathbb{E}\left[\exp{-\langle
k\bar{\phi},J_t\rangle}\right]\mid_{k=0^+}\\& = &
\left\langle\frac{d}{dk}V_t^0(k\bar{\phi}),F\times\Lambda\right\rangle\exp\left({-\left\langle
V_t^0(k\bar{\phi}),F\times\Lambda\right\rangle}\right)\bigg|_{k=0^+}.
\end{eqnarray*}
Thus, putting
$\dot{V}_t\bar{\phi}:=\frac{d}{dk}V_t^0(k\bar{\phi})\mid_{k=0^+}$
and recalling that $V_t^0(0\bar{\phi})=0$, we obtain
\begin{center}
\begin{tabular}{lll}
$\frac{\partial}{\partial t}\dot{V}_t\bar{\phi}(u,x)$ & = & $\mathcal{L}\dot{V}_t\bar{\phi}(u,x)+\bar{\phi}(u,x)$\\
$\dot{V}_0\bar{\phi}(u,x)$ & = & 0,
\end{tabular}
\end{center}
or
\begin{equation}\label{asterisco4}
\dot{V}_t\bar{\phi}(u,x)=\int_0^t\tilde{T}_{t-s}\bar{\phi}(u,x)ds=\int_0^t\S_{t-s}\phi(x)ds.
\end{equation}
Consequently, using that $\Lambda$ is invariant for the
$\alpha$-stable semigroup,
\begin{equation*}
\mathbb{E}\langle\bar{\phi},J_t\rangle=
\langle\dot{V}_t\bar{\phi},F\times\Lambda\rangle=
\left\langle\int_0^t\S_{t-s}{\phi}
ds,\Lambda\right\rangle=\left\langle{\phi},\Lambda\right\rangle t.
\end{equation*}
This proves (\ref{asterisco3}).
The proof of  (\ref{CotaVar}) goes as follows. Define $\bar{\phi}$ as before.
 Differentiating $V_t(k\bar{\phi})$ with respect to $k$ and using
equation (\ref{TOC}), we obtain
\begin{eqnarray*}
\frac{\partial^2}{\partial t\partial k}V_t(k\bar{\phi})(u,x)& = &
\L\frac{\partial}{\partial
k}V_t(k\bar{\phi})(u,x)+\bar{\phi}(u,x)(1-V_t(k\bar{\phi})(0,x))-k\bar{\phi}(u,x)\frac{\partial}{\partial
k}V_t(k\bar{\phi})(u,x)\\&  &
-\lambda(u)\left[-\Phi^\prime(1-V_t(k\bar{\phi}(0,x))\frac{\partial}{\partial
k}V_t(k\bar{\phi})(0,x)+\frac{\partial}{\partial
k}V_t(k\bar{\phi})(0,x)\right],
\end{eqnarray*}
and
\begin{eqnarray*}
\frac{\partial^3}{\partial t\partial k^2}V_t(k\bar{\phi})(u,x)& = &
\L\frac{\partial^2}{\partial
k^2}V_t(k\bar{\phi})(u,x)-2\bar{\phi}(u,x)\frac{\partial }{\partial
k}V_t(k\bar{\phi})(u,x)-k\bar{\phi}(u,x)\frac{\partial^2}{\partial
k^2}V_t(k\bar{\phi})(u,x)\\&  &
-\lambda(u)\left[\Phi^{\prime\prime}\left(1-V_t(k\bar{\phi})(0,x)\right)\left(\frac{\partial}{\partial
k}V_t(k\bar{\phi})(u,x)\right)^2\right.\\&  & \left.
-\Phi^\prime(1-V_t(k\bar{\phi})(0,x))\frac{\partial^2}{\partial
k^2}V_t(k\bar{\phi})(0,x)+\frac{\partial^2}{\partial
k^2}V_t(k\bar{\phi})(0,x)\right].
\end{eqnarray*}
Letting $\ddot{V}_t\bar{\phi}=\frac{\partial^2}{\partial
k^2}V_t(k\bar{\phi})|_{k=0^+}$ we get that
\begin{equation}\label{DOS}
\frac{\partial}{\partial t}\ddot{V}_t\phi(u,x) =
\L\ddot{V}_t\bar{\phi}(u,x) -
\lambda(u)\Phi^{\prime\prime}(1)\left(\dot{V}_t\bar{\phi}(0,x)\right)^2
- 2\bar{\phi}(u,x)\dot{V}_t\bar{\phi}(u,x).
\end{equation}
From (\ref{asterisco4}) and (\ref{DOS}) we obtain
\begin{equation*}
\ddot{V}_t\bar{\phi}(u,x) =
\int_0^t\tilde{T}_s\left[-\lambda(u)\left(\int_0^s\tilde{T}_r\bar{\phi}(u,x)dr\right)^2
- 2\bar{\phi}(u,x)\int_0^s\tilde{T}_r\bar{\phi}(u,x)dr \right]ds.
\end{equation*}
Note that $\mbox{Var}\langle\phi,J_t\rangle = -\langle
\ddot{V}_t\bar{\phi}(*,\bullet),F\times\Lambda\rangle$. Therefore,
\begin{eqnarray}\label{aux}
\nonumber\mbox{Var}\langle\bar{\phi},J_t\rangle& = & \left\langle
\int_0^t\tilde{T}_s\left[\lambda(*)\left(\int_0^s\tilde{T}_r\bar{\phi}
(*,\bullet)dr\right)^2 +
2\bar{\phi}(*,\bullet)\int_0^s\tilde{T}_r\bar{\phi}(*,\bullet)\,dr
\right]\,ds,F\times\Lambda\right\rangle\\ \nonumber& = &
\int_0^t\left\langle \lambda(*)\left(\int_0^s\tilde{T}_r\bar{\phi}
(*,\bullet)dr\right)^2,F\times\Lambda\right\rangle ds\\\nonumber &
& + 2\int_0^t\left\langle
\bar{\phi}(*,\bullet)\int_0^s\tilde{T}_r\bar{\phi}
(*,\bullet)dr,F\times\Lambda\right\rangle ds\\&=:& (A)+(B).
\end{eqnarray}
Notice that, under the choice of $\bar{\phi}$,
$\tilde{T}_t\bar{\phi}(u,x)=\mathcal{S}_t\phi(x)$ for all $t\geq0$,
and that $\langle\lambda,F\rangle<\infty$. In fact, using that
$\lambda(u)\sim u^{-1}$ and $f(u)\sim u^{-\gamma-1}$, we get that,
for $A>0$ sufficiently large,
\begin{eqnarray*}
\langle\lambda,F\rangle & =
&\int_0^\infty\lambda(u)f(u)du\\&=&\int_0^A\lambda(u)f(u)du+\int_A^\infty\lambda(u)f(u)du\\&\sim&\int_0^A\lambda(u)f(u)du+\int_A^\infty
u^{-1}u^{-\gamma-1}du\\&<&\infty.
\end{eqnarray*}
Now,
$$
(A)  =  \int_0^t\left\langle
\lambda(*)\left(\int_0^s\S_r\phi(\bullet)\,dr
\right)^2,F\times\Lambda\right\rangle\, ds =
\langle\lambda,F\rangle\int_0^t\left \langle\left(\int_0^s\S_r\phi
dr\right)^2,\Lambda\right\rangle\, ds.
$$
Also, it can be shown that
\begin{equation*}
\int_0^t\left\langle\left(\int_0^s\S_r\phi
dr\right)^2,\Lambda\right\rangle ds\leq \mbox{Const}(\phi_2)(t +
t^{3-d/\alpha}),
\end{equation*}
and consequently,
\begin{equation*}
(A) \leq \langle\lambda,F\rangle\mbox{Const}(\phi_2)(t +
t^{3-d/\alpha}).
\end{equation*}
Similarly, for the second term in (\ref{aux})
$$
(B) =  2\int_0^t\left\langle
\bar{\phi}(*,\bullet)\int_0^s\tilde{T}_r\bar{\phi}(*,\bullet)\,
dr,F\times\Lambda\right\rangle \,ds=
2\int_0^t\left\langle\int_0^s\S_r\phi dr,\Lambda\right\rangle\, ds,
$$
where
\begin{equation*}
\int_0^t\left\langle\int_0^s\S_r\phi\,dr,\Lambda\right\rangle
ds\leq\mbox{Const}(\phi)(t + t^{2-d/\alpha}),
\end{equation*}
hence,
\begin{equation*}
(B) \leq 2\mbox{Const}(\phi)(t + t^{2-d/\alpha}).
\end{equation*}
Finally, combining  the bounds for (A) and (B) we obtain the
result.\hfill$\Box$

\section{Proof of Theorem \ref{LLNFinitemean}}
Suppose that $F$ is a non-arithmetic distribution function supported
on the non-negative real line  and having finite mean $\mu>0$, and let
$d>\alpha$.
As in the
proof of Lemma \ref{cota}, we have that
\begin{equation}
P\left(|T^{-1}\langle\varphi,J_T(1)\rangle-\langle\varphi,\Lambda\rangle|>\epsilon\right)
\leq\frac{2}{\epsilon^2}\int_0^1dv\int_0^v\mbox{Cov}\left(\langle\varphi,X_{uT}\rangle,\langle\varphi,X_{vT}\rangle\right)du.
\end{equation}
Therefore, due to Proposition \ref{COVp},
\begin{equation}\label{aux1}
P\left(|T^{-1}\langle\varphi,J_T(1)\rangle-\langle\varphi,\Lambda\rangle|>\epsilon\right)\leq
(I)+(II),
\end{equation}
where
\begin{equation*}
(I):=\frac{2}{\epsilon^2}\int_0^1dv\int_0^vdu\langle\varphi\S_{T(v-u)}\varphi,\Lambda\rangle,
\end{equation*}
and
\begin{equation*}
(II):=\frac{2}{\epsilon^2}\int_0^1\int_0^v\int_0^{Tu}\langle(\S_{Tu-r}\varphi)(\S_{Tv-r}\varphi),\Lambda\rangle
dU(r)du\,dv.
\end{equation*}
We recall the bound (\ref{cotaI}) for $(I)$. It remains to
upper-bound $(II)$. Performing the change of variables $h=r/T$ in
$(II)$, and using the {elementary renewal theorem} (see e.g.
\cite{KT}, p. 188), we have that, for $T$ large enough,
\begin{eqnarray*}
(II)&=&\frac{2}{\epsilon^2}\int_0^1\int_0^v\int_0^{u}
\langle(\S_{T(u-h)}\varphi)(\S_{T(v-h)}\varphi),\Lambda\rangle
d\left[\frac{U(Th)}{Th}Th\right]du\,dv\\&\sim&\frac{2T}{\epsilon^2\mu}\int_0^1\int_0^v\int_0^{u}
\langle(\S_{T(u-h)}\varphi)(\S_{T(v-h)}\varphi),\Lambda\rangle
dh\,du\,dv\\&=&\frac{2T}{\epsilon^2\mu}\int_0^1\int_0^v\int_h^v\langle(\S_{T(u-h)}
\varphi)(\S_{T(v-h)}\varphi),\Lambda\rangle
du\,dh\,dv.
\end{eqnarray*}
After performing several changes of variables one can see that, for all
$T$ large enough,
\begin{eqnarray*}
(II)&\sim&\frac{2}{\epsilon^2\mu
T^2}\int_0^T\int_{\Rd}\int_0^v\int_0^t(\S_s\varphi)(x)(\S_t\varphi)(x)
ds\,dt\,dx\,dv \\&\leq& \frac{2}{\epsilon^2\mu
T^2}\int_0^T\int_{\Rd}\int_0^v\int_0^v(\S_s\varphi)(x)(\S_t\varphi)(x)ds\,dt
\,dx\,dv\\&=&\frac{2}{\epsilon^2\mu
T^2}\int_0^T\int_{\Rd}\int_{\Rd}\varphi(y)\varphi(z)\int_0^v\int_0^vp_{t+s}(y-z)ds\,dt\,dy\,dz\,dv.
\end{eqnarray*}
On the other hand,   one can show, as in \cite{SS}, that
\begin{equation*}
\int_0^v\int_0^vp_{t+s}(y-z)ds\,dt\leq
c\left(|y-z|^{2\alpha-d}+v^{2-d/\alpha}\right)
\end{equation*}
for some constant $c>0$. Hence, for any fixed $A>0$, and all $T$
large enough,
\begin{eqnarray*}
(II)&\leq&\frac{2}{\epsilon^2\mu
T^2}\int_0^A\int_{\Rd}\int_{\Rd}\varphi(y)\varphi(z)\int_0^v\int_0^vp_{t+s}(y-z)ds\,dt\,dy\,dz\,dv\\&
&+\frac{2}{\epsilon^2\mu
T^2}c\int_{\Rd}\int_{\Rd}\varphi(y)\varphi(z)|y-z|^{2\alpha-d}dy\,dz(T-A)\\&
&+\frac{2}{\epsilon^2\mu
T^2}c\int_{\Rd}\int_{\Rd}\varphi(y)\varphi(z)dy\,dz\frac{(T^{3-d/\alpha}-A^{3-d/\alpha})}{3-d/\alpha}\\&\leq&\frac{2}{\epsilon^2\mu
T^2}\int_0^A\int_{\Rd}\int_{\Rd}\varphi(y)\varphi(z)\int_0^v\int_0^vp_{t+s}(y-z)ds\,dt\,dy\,dz\,dv\\&
&+\frac{2}{\epsilon^2\mu T
}c\int_{\Rd}\int_{\Rd}\varphi(y)\varphi(z)|y-z|^{2\alpha-d}dy\,dz\\&
&+\frac{2}{\epsilon^2\mu
}c\int_{\Rd}\int_{\Rd}\varphi(y)\varphi(z)dy\,dz\frac{(T^{1-d/\alpha}-A^{3-d/\alpha}T^{-2})}{3-d/\alpha}.
\end{eqnarray*}
The proof concludes with an application of   Borel-Cantelli's Lemma,
using that $d/\alpha>1$, (\ref{aux1}), and the bounds for (I) and
(II).

\section{Proof of Theorem \ref{subcritical}}
In this section we assume that $d<\alpha\gamma$, and that $F$ is a
distribution function satisfying (\ref{tail1}). Arguing similarly as
at the end of the proof of Theorem \ref{lowdimensionsT}, Lemma
\ref{BCL} below yields the theorem.

\begin{lemma} \label{BCL}Let $A\subset\Rd$
be a ball. Then, for all $\epsilon>0$, and for all $t$
sufficiently large,
\begin{equation}\label{BCA}
 P\left(t^{-1}\int_0^t1_{\{X_s(A)>0\}}ds>\epsilon\right)\leq(1-e^{-\epsilon})^{-1}
 \left(ct^{-(d/\alpha+\gamma)/2+(1+\delta)d/\alpha}+c_1t^{-1}\right),
\end{equation}
for some positive constants $c$ and $c_1$. \end{lemma}
 \demostracion Notice that, by Markov's inequality,
 \begin{equation}\label{BC}
 P\left(t^{-1}\int_0^t1_{\{X_s(A)>0\}}ds>\epsilon\right)\leq(1-e^{-\epsilon})^{-1}\E\left[1-\exp\left\{t^{-1}
 \int_0^t1_{\{X_{s}(A)>0\}}ds\right\}\right].
 \end{equation}
 Moreover, since
the initial population is Poissonian,
 \begin{equation}\label{Exp}
 \E\left[e^{-t^{-1}\int_0^t1_{\{X_s(A)>0\}}ds}\right]=
 \E\left[e^{-\int_0^11_{\{X_{st}(A)>0\}}ds}\right]=\exp\left\{-\int_{\Rd}\Ex\left[1-e^{-\int_0^11_{\{Z_{st}(A)>0\}}ds}\right]dx\right\}.
 \end{equation}
 Now, since $1-e^{-x}\leq x$, for all $x\geq0$, we have that
 \begin{equation*}
 1-e^{-\int_0^11_{\{Z_{st}(A)>0\}}ds}\leq \int_0^11_{\{Z_{st}(A)>0\}}ds.
 \end{equation*}
Therefore,
 \begin{equation}\label{exp} \Ex\left[1-e^{-\int_0^11_{\{Z_{st}(A)>0\}}ds}\right]\leq
\Ex\int_0^11_{\{Z_{st}(A)>0\}}\,ds=\int_0^1P_x\left(Z_{st}(A)>0\right)\,ds.
 \end{equation}
 Due to (\ref{Exp}) and (\ref{exp}), inequality (\ref{BC}) can be written as
\begin{equation}\label{BC1}
P\left(t^{-1}\int_0^t1_{\{X_{s}(A)>0\}}ds>\epsilon\right)\leq(1-e^{-\epsilon})^{-1}\left[1-\exp
\left\{-\int_0^1\int_{\Rd}P_x\left(Z_{st}(A)>0\right)dx\,ds\right\}\right],
\end{equation}
where
\begin{equation}\label{majorized}
\int_0^1\int_{\Rd}P_x\left(Z_{st}(A)>0\right)dx\,ds=\int_0^1\int_{D(st,\delta)}P_x
\left(Z_{st}(A)>0\right)dx\,ds+\int_0^1\int_{D(st,\delta)^c}P_x\left(Z_{st}(A)>0\right)dx\,ds
\end{equation}
with
\begin{equation*}
D(t,\delta):=\{x\in\Rd: |x|\leq
t^{(1+\delta)/\alpha}\},\;\;\delta>0.
\end{equation*}
Using the inequality
\begin{equation*}
\int_{D(st,\delta)}P_x\left(Z_{st}(A)>0\right) dx\leq
K(st)^{-(d/\alpha+\gamma)/2+(1+\delta)d/\alpha},
\end{equation*}
which holds for some positive constant $K$ (see  Lemma 5 in \cite{VW}), we
deduce that
\begin{eqnarray*}
\int_0^1\int_{D(st,\delta)}P_x\left(Z_{st}(A)>0\right) dx\,ds&\leq&
K\int_0^1(st)^{-(d/\alpha+\gamma)/2+(1+\delta)d/\alpha}
ds\\&=&ct^{-(d/\alpha+\gamma)/2+(1+\delta)d/\alpha},
\end{eqnarray*}
where we used that $1{-(d/\alpha+\gamma)/2+(1+\delta)d/\alpha}>0$. On the other
hand, following closely the proof of Lemma 5 in \cite{VW} one can
see that, for sufficiently large $t$,
$$
\int_0^1\int_{\Rd\backslash D(st,\delta)}P_x\left(Z_{st}(A)>0\right)dx\,ds\leq
c_1\int_0^1P\left(\|B^0_1\|\geq\frac{1}{2}(st)^{\delta/\alpha}\right)\,ds
\leq
c_1t^{-1}.
$$
In this way, (\ref{majorized}) yields the inequality
\begin{equation}\label{majorized1}
\int_0^1\int_{\Rd}P_x\left(Z_{st}(A)>0\right)dx\,ds\leq
ct^{-(d/\alpha+\gamma)/2+(1+\delta)d/\alpha}+c_1t^{-1},
\end{equation}
which is valid for all $t$ large enough, and renders (\ref{BCA}).
Notice that ${-(d/\alpha+\gamma)/2+(1+\delta)d/\alpha}<0$ for
sufficiently small $\delta$.\hfill$\Box$\medskip

{\bf\noindent Acknowledgement}
The authors express their gratitude to an anonymous referee for her/his meti-\\culous revision of the paper, and for pointing out a mistake in the proof of an earlier version of Theorem \ref{subcritical}.

\end{document}